\documentclass{article}

\usepackage[english]{babel}

\usepackage[letterpaper,top=2cm,bottom=2cm,left=3cm,right=3cm,marginparwidth=1.75cm]{geometry}

\usepackage{amsmath}
\usepackage{graphicx}
\usepackage[colorlinks=true, allcolors=blue]{hyperref}

\makeatletter
\patchcmd{\@verbatim}
  {\verbatim@font}
  {\verbatim@font\small}
  {}{}
\makeatother

\usepackage{amsthm}
\theoremstyle{definition}
\newtheorem{conjecture}{Conjecture}

\title{Irregular Stanley sequences plausibly do not have growth $\Theta(n^2/\log n)$}
\author{Nat Sothanaphan\footnote{natsothanaphan@gmail.com}}
\date{December 13, 2025}

\begin{document}
\maketitle

\begin{abstract}
Stanley sequences starting from the set $\{0, n\}$ where $n$ is a positive integer have long been conjectured to be divided into two types: the ``regular'' type where the growth rate is $\Theta(n^{\log_2(3)})$, and the ``irregular'' type where the growth rate is thought to be $\Theta(n^2/\log n)$. A paradigmatic case of a candidate irregular type is $n=4$, although to date no value of $n$ has been proven to have such a growth rate. Here, we provide strong numerical evidence \emph{against} this conjectured growth rate for $n=4$. Specifically, for $n=4$, it seems plausible that the upper bound is $O(n^2/\log n)$ but that the lower bound is in fact $\Omega(n^{2-\delta})$ for some $\delta > 0$. This appears to be because the sequence is not totally ``random'' as has been assumed. Limitations of the numerical method here is discussed.
\end{abstract}

\section{Introduction}

Stanley sequences are introduced by Odlyzko and Stanley \cite{odlyzko-stanley-1978-greedy} in their 1978 paper. There, they gave a heuristic probabilistic argument that if the sequence behaves ``randomly,'' then the growth rate should be $\Theta(n^2/\log n)$. They also noted that the sequence for an irregular starting value behaves ``erratically.'' This influences the perspectives of later authors; indeed, regular sequences are much more often studied and the irregular ones only cited as being erratic or chaotic.

This short note provides numerical evidence for the following conjecture.

\begin{conjecture}
The Stanley sequence starting from $\{0, 4\}$ (an ``irregular'' case) does not have growth rate $\Theta(n^2/\log n)$ as previously assumed. Instead, the upper bound $O(n^2/\log n)$ does hold, but the lower bound is approximately $\Omega(n^{2-\delta})$ for some $\delta > 0$ (meaning that the infimum of $\delta$ that makes this hold is strictly positive).
\end{conjecture}

The author thanks ChatGPT for numerical investigation.

\section{Previous literature}

Growth rates of Stanley sequences have been studied by Moy \cite{moy-2011-counting-stanley}, Dai and Chen \cite{dai-chen-2013-counting-stanley}, and Rolnick and Venkataramana \cite{rolnick-venkataramana-2015-growth-stanley}. Erdos had an interest in the problem. Indeed, the growth rate of Stanley sequences starting from $\{0, n\}$ is listed as Problem 271 in the Erdos Problems website \cite{bloom-erdos-271}, and this remains an open problem.

Granted, limitations on computational power back then could have limited experiments. However, we note that Lindhurst \cite{lindhurst-1990-greedy-thesis} performed computation on the case $n=4$, and that this case also appears as sequence A005487 in the OEIS \cite{oeis-A005487}. Hence, one must assume that there is extensive data supporting the case that $n=4$ gives the growth rate of $\Theta(n^2/\log n)$. Imagine my surprise when this turns out to not be the case!

\section{Numerical exploration}

It is probably useful to describe the algorithm, so that other people can also extend the experiments here.

\subsection{Algorithm}

A Stanley sequence is just a sequence constructed in a greedy manner to contain no $3$-term arithmetic progression and such that the next element must be greater than the previous one. This means no $a+c=2b$. So to check whether $c$ should be included, we check $2b - a$ where $b > a$.

In practice, it's faster to check whether $2b-c$ is in the sequence thus far (``seen''). If we use a set which has $O(1)$ lookup via hashing, this can be fast. Also, an important optimization is that we check by enumerating $b$ in \emph{decreasing} rather than increasing order: this becomes clear once we realize that if we have, say, $n$ and $n+1$, we should quickly rule out $n+2$ by going downward rather than upward from $0$!

\subsection{Potential fractal-like structure}

At first the author (and ChatGPT) expected a ``random'' behavior, as discussed elsewhere. However, something different occurs.

\begin{figure}
\centering
\includegraphics[width=250px]{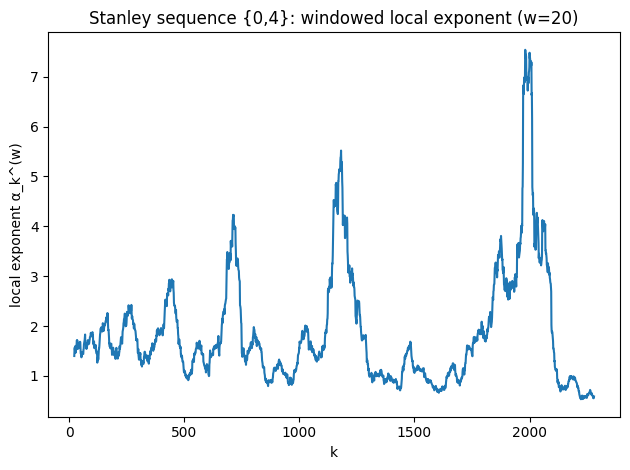}
\caption{\label{fig:windowed-exponent}Windowed local exponent reveals potential self-similar structure.}
\end{figure}

We're looking at the ``windowed exponent'' defined by
$$ \alpha_k^{(w)} = \frac{\log a_{k+w} - \log a_{k-w}}{\log (k+w) - \log (k-w)}. $$
This is a form of local computation of the exponent if we expect $a_k \approx k^{\alpha}$ where $\alpha$ should somehow cluster around $2$ for $k$ large. Instead as observed in Figure \ref{fig:windowed-exponent}, there is increasingly divergent behavior along with a potential self-similarity or fractal-like structure.

This is reminiscent of the binary-ternary explicit representation of the Stanley sequence starting from $\{0, 1\}$ (a regular case), where there are jumps at multiple, exponentially spaced scales. It hints that the regular and irregular cases may in fact not be so different as previously thought.

However, theoretical analysis remains difficult and has not been attempted here.

\subsection{Peaks and troughs}

Now, recall at this point Odlyzko and Stanley's argument that the growth rate of $\Theta(n^2/\log n)$ would follow if the sequence behaves sufficiently randomly. But we have found evidence that it doesn't behave randomly! Hence the immediate question: do we still have this conjectured growth rate?

What follows is a rather nuanced attempt at numerical analysis. Notice that there should be exponentially-spaced ``peaks and troughs'' on the graph. So we should look at these subsequences. Also, we need to look at the plot of $\log a_k/\log k$ directly rather than the windowed exponent.

It's not that straightforward to locate peaks/troughs programmatically. But here's one way. First we smooth out by convolution. Then we use scipy's \verb|find_peaks| method. It may not be perfect but should produce candidates to be further selected or added by hand. In any case, we now have peaks and troughs.

\begin{figure}
\centering
\includegraphics[width=300px]{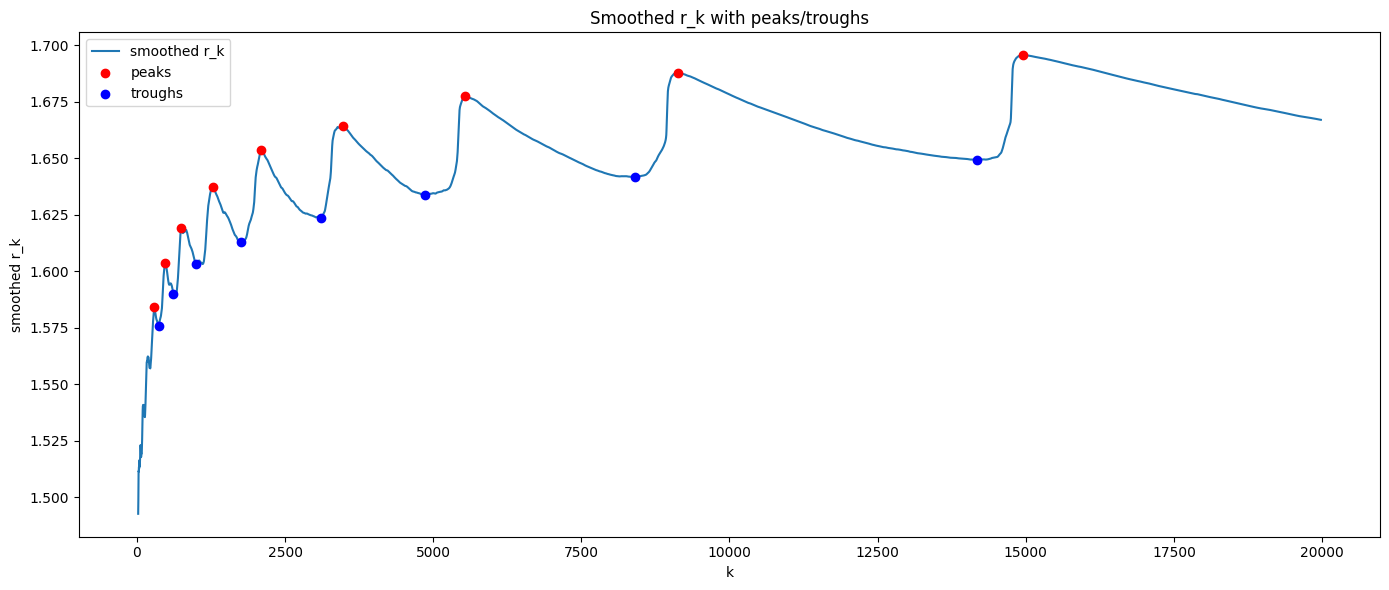}
\caption{\label{fig:peaks-and-troughs}Peaks and troughs of $\log a_k/\log k$.}
\end{figure}

See Figure \ref{fig:peaks-and-troughs}.

\subsection{Growth rates}

It now remains to fit the peak/trough sequences to see if they behave like $\Theta(n^2/\log n)$. Visually, I initially tried to overlay $Cn^2/\log n$ directly, $C$ being adjustable constant. That seems to fit the ``peaks'' well, but the ``troughs'' appear too flat to fit. This leads to the conjecture that $\Theta(n^2/\log n)$ is indeed satisfied as an \emph{upper bound} but not a lower bound.

Nevertheless, this is not strong evidence as it depends on possible human bias. Hence we employ the following method. We try to fit a linear regression by finding best $A,B,C$ such that
$$ \log a_k \approx A + B \frac{\log \log k}{\log k} + C \frac{1}{\log k}. $$
So for peaks, we expect $A \approx 2$ and for troughs, perhaps $A < 2$.

While it seems like a good idea, I discovered that there is simply too much ``noise'' in the values to estimate the coefficients robustly. We have $9$ peaks and $8$ troughs here (already selected so the $k$ is not too small to skew the curve too much). To test for robustness, I regressed three times: the whole sequence, when deleting the last element, and when deleting the first element. For the peaks, I get $A \in \{1.947, 1.984, 1.900\}$ and for the troughs, I get $A \in \{1.826, 1816, 1.893\}$ (rounded to three decimal places). So it seems like there is a consistent gap $A_{\text{trough}} < A_{\text{peak}}$. However, we see that there is too much noise to conclude this.

The next attempt was to fix $A=2$ and regress to find only $B,C$. Here the results are more robust: for peaks $B \in \{-1.050, -1.031, -1.054\}$; for troughs $B \in \{-1.730, -1.698, -1.805\}$. Note how $\Theta(n^2/\log n)$ implies $A=2, B=-1$. This is extremely consistent for peaks, but not consistent for troughs. Thus, we have the first moderately strong evidence that the troughs do not satisfy $\Theta(n^2/\log n)$ while the peaks do. (Note on interpretation: we believe that $A_{\text{trough}} < 2$, but the (statistical) argument here is that \emph{if} $A_{\text{trough}} = 2$, then it is unlikely that $B_{\text{trough}} = -1$. It is sufficient to go against this specific combination.)

Full coefficients are in the Appendix.

\section{Evidence against $\Theta(n^2/\log n)$}

\subsection{Direct plot}

\begin{figure}
\centering
\includegraphics[width=350px]{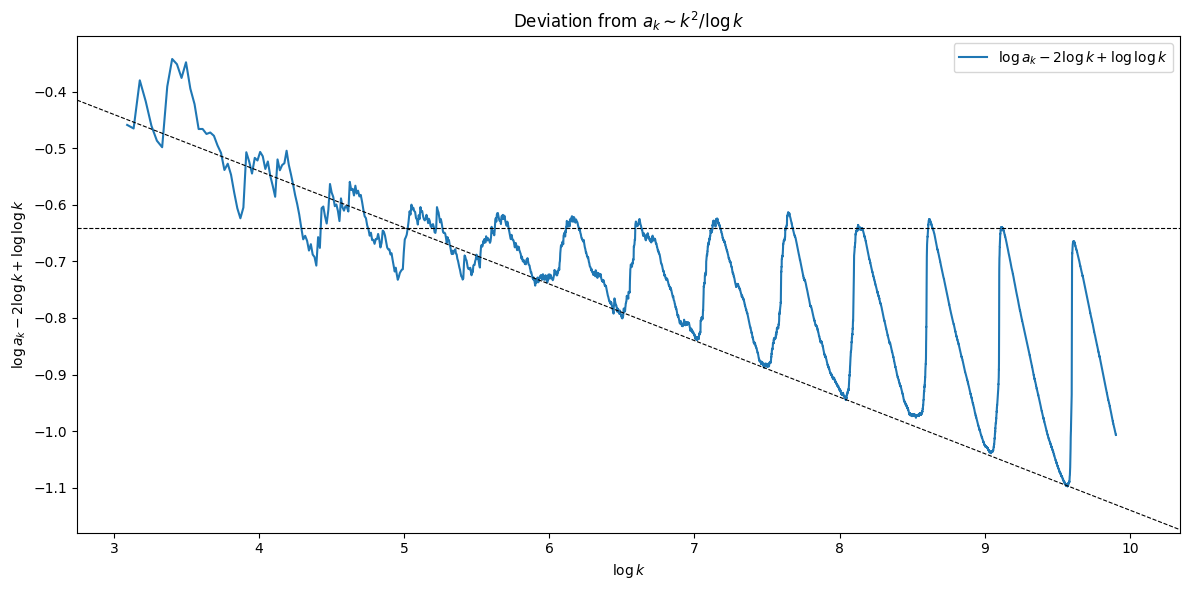}
\caption{\label{fig:direct-plot}Direct plot of $\log k$ vs. $\log a_k - 2 \log k + \log \log k$.}
\end{figure}

The exploration thus far does provide evidence against $\Theta(n^2/\log n)$ in the lower bound; however, it is relatively convoluted. Of course, once we know what we're looking for, we can set up a direct plot:

We plot
$$ \log k \text{ vs. } \log a_k - 2 \log k + \log \log k $$
directly. The conjectured growth rate of $\Theta(n^2/\log n)$ implies that the graph must be eventually bounded. The result is as in Figure \ref{fig:direct-plot}.

Some information about this plot:
\begin{itemize}
\item The first 20 values are skipped to make it not fluctuate too much.
\item We have data up to $k = 20\text{,}000$. Large data is difficult obtain because of long runtime (currently about 16 minutes in Colab for this $k$).
\item The horizontal ``upper bound'' line is $y = -0.64$.
\item The slanted ``lower bound'' line is $y = -0.1 x - 0.14$.
\end{itemize}

Here's what we can conclude: the ``bounded above'' story is well supported by data, indicated with the manually drawn horizontal line.

Nevertheless, the ``bounded below'' story is not well supported. In particular, there is regularity in the plot where peaks and troughs seem to appear at equally spaced $\log k$ points (which would be exponentially spaced $k$ points). And the troughs seem be going down linearly in the plotted coordinates with high regularity. If the slope of $-0.1$ is a good approximation, then the correct lower bound could be approximately $n^{1.9}$.

\subsection{Plausible objections}

Of course, we should take a skeptical mindset.

\emph{Objection 1}. This region is too small for asymptotics to kick in.

\emph{Reply 1}. In our investigation of growth rates, where we fit regression with $A=2$, we have obtained $B \approx -1$ for peaks but $B < -1$ for troughs. Since the value for peaks aligns well with the theoretical bound, this is evidence that we're already within the region where asymptotics make sense, at least in its ability to distinguish the growth rates of the peak/trough subsequences.

\emph{Objection 2}. Perhaps we should trust the $\Theta(n^2/\log n)$ heuristic more because it accurately predicts the peaks; maybe the troughs are indeed bounded but for reasons we haven't considered.

\emph{Reply 2}. Due to the non-random fractal-like structure observed empirically in Figure \ref{fig:windowed-exponent}, it seems more likely that the heuristic used to obtain the $\Theta(n^2/\log n)$ is inaccurate instead. Indeed, the observed peaks and troughs really do occur with high regularity.

So, taken together the direct plot and the exploration previously done, we have obtained a strong case against the $\Theta(n^2/\log n)$ growth rate even with limited data.

\subsection{Limitations of the method}

For one, sequences with no $3$-term arithmetic progression are difficult to analyze theoretically. This suggests numerical experimentations as we have done.

But even with this, there is great limitation. Our algorithm is somewhere between $O(a_n)$ and $O(n a_n)$ in runtime (depending on whether an element is quickly ruled out); so probably somewhere between $O(n^2)$ and $O(n^3)$. However, peaks and troughs grow exponentially. So we can realistically obtain only a few more peaks and troughs even if we throw in more compute.

As we have seen, within this regime, there is enough noise in the peak/trough data that prevents robust regression. This means that we're not likely to be able to estimate the correct bound from the numerical data!

Based on this limitation, I have thus chosen to only conclude that the \emph{likely} upper bound appears to be $O(n^2/\log n)$ as theoretically predicted but that the lower bound is around $\Omega(n^{2-\delta})$ for some $\delta > 0$ (possibly around $0.1$). The evidence presented is the combination of many analyses done rather than just one, as described in the previous subsection, to strengthen the claim from limited data.

\section{Conclusion}

In any case, that the bound for irregular Stanley sequences starting from $\{0, n\}$ may not be as predicted and that there may be structures in them surely catch me off guard. I hope these investigations inspire further works on the irregular case.

\bibliographystyle{alpha}
\bibliography{references}

\section*{Appendix: code}

Here is the Python code for reference. It can be run in Colab environment.

\subsection*{Sequence generation}

\begin{verbatim}
import numpy as np
import math
import matplotlib.pyplot as plt

# ---------- Stanley sequence {0, n0} ----------

def stanley_sequence_two_seed(n0, max_len):
    """Greedy 3-AP-free Stanley sequence starting from {0, n0}."""
    seq = [0, n0]
    seen = set(seq)

    def is_valid(c):
        # adding c must not create a 3-term AP a < b < c
        for i in range(len(seq)-1):
            b = seq[-i-1] # it's faster to check in decreasing order
            a = 2*b - c
            if a in seen:
                return False
            if a < 0:
                return True
        return True

    while len(seq) < max_len:
        candidate = seq[-1] + 1
        while True:
            if is_valid(candidate):
                if len(seq) % 1000 == 0:
                    print("we're at length: ", len(seq))
                seq.append(candidate)
                seen.add(candidate)
                break
            candidate += 1
    return np.array(seq, dtype=int)

# ---------- parameters ----------

N = 20000
# N=20000 took about 16 mins in Colab environment
n0 = 4

seq = stanley_sequence_two_seed(n0, N)
\end{verbatim}

\subsection*{Windowed exponent}

\begin{verbatim}
# Windowed local exponent alpha_k^(w)
w = 20
ks_loc = list(range(w+2, N-w+1))  # ensure indices stay in range and avoid a_1=0
alphas = []
for k in ks_loc:
    j_minus = k - w
    j_plus = k + w
    num = math.log(seq[j_plus-1]) - math.log(seq[j_minus-1])
    den = math.log(j_plus) - math.log(j_minus)
    alphas.append(num / den)

plt.figure()
plt.plot(ks_loc, alphas)
plt.xlabel("k")
plt.ylabel("local exponent alpha_k^(w)")
plt.title(f"Stanley sequence {{0,4}}: windowed local exponent (w={w})")
plt.tight_layout()
plt.show()
\end{verbatim}

\subsection*{Finding peaks and troughs}

Part 1:

\begin{verbatim}
# r_k = log(a_k)/log(k) for k >= 2
ks = np.arange(2, N+1)
a_vals = seq[1:]                 # a_2,...,a_N
r_vals = np.log(a_vals) / np.log(ks)
\end{verbatim}

Part 2:

\begin{verbatim}
# ---------- smoothing r_k ----------

L = 25                # smoothing window length; adjust by eye
kernel = np.ones(L) / L
r_smooth = np.convolve(r_vals, kernel, mode="same")

plt.figure(figsize=(10,4))
plt.plot(ks[20:-20], r_vals[20:-20], alpha=0.3, label="raw r_k")
plt.plot(ks[20:-20], r_smooth[20:-20], label=f"smoothed r_k (L={L})")
plt.xlabel("k")
plt.ylabel("r_k = log(a_k)/log k")
plt.title("Exponent estimate r_k and smoothed version")
plt.legend()
plt.tight_layout()
plt.show()
\end{verbatim}

Part 3:

\begin{verbatim}
from scipy.signal import find_peaks

# ---------- candidate peaks/troughs on smoothed r_k ----------

# tune these:
min_dist = 50    # minimum separation in k between picked peaks
prom_peak = 0.005 # required prominence for peaks (radians-ish; small here)
prom_trough = 0.003

# candidate peaks: high points of r_smooth
peak_idx, peak_props = find_peaks(r_smooth,
                                  distance=min_dist,
                                  prominence=prom_peak)

# candidate troughs: high points of -r_smooth
trough_idx, trough_props = find_peaks(-r_smooth,
                                      distance=min_dist,
                                      prominence=prom_trough)

k_peaks_cand = ks[peak_idx]
r_peaks_cand = r_smooth[peak_idx]

k_troughs_cand = ks[trough_idx]
r_troughs_cand = r_smooth[trough_idx]

print("candidate peak ks:", k_peaks_cand)
print("candidate trough ks:", k_troughs_cand)

plt.figure(figsize=(10,4))
plt.plot(ks[20:-20], r_smooth[20:-20], label="smoothed r_k")
plt.scatter(k_peaks_cand, r_peaks_cand, color="red", zorder=3, label="candidate peaks")
plt.scatter(k_troughs_cand, r_troughs_cand, color="blue", zorder=3, label="candidate troughs")
plt.xlabel("k")
plt.ylabel("smoothed r_k")
plt.legend()
plt.title("Smoothed r_k with candidate peaks/troughs")
plt.tight_layout()
plt.show()

r_peaks_cand, r_troughs_cand
\end{verbatim}

Part 4:

\begin{verbatim}
k_peaks = [293, 480, 750, 1285, 2100, 3486, 5538, 9131, 14957] # pick manually
k_troughs = [365, 618, 1001, 1765, 3107, 4854, 8410, 14179] # pick manually
r_peaks = [r_smooth[k-2] for k in k_peaks]
r_troughs = [r_smooth[k-2] for k in k_troughs]

plt.figure(figsize=(14,6))
plt.plot(ks[20:-20], r_smooth[20:-20], label="smoothed r_k")
plt.scatter(k_peaks, r_peaks, color="red", zorder=3, label="peaks")
plt.scatter(k_troughs, r_troughs, color="blue", zorder=3, label="troughs")
plt.xlabel("k")
plt.ylabel("smoothed r_k")
plt.legend()
plt.title("Smoothed r_k with peaks/troughs")
plt.tight_layout()
plt.show()

r_peaks, r_troughs
\end{verbatim}

\subsection*{Regression}

Part 1:

\begin{verbatim}
# this is just a rename here
peak_k = np.array(k_peaks)
trough_k = np.array(k_troughs)

# compute r_k at those exact indices using the *unsmoothed* r_vals
def r_at(k):
    return math.log(seq4[k-1]) / math.log(k)

peak_r = np.array([r_at(k) for k in peak_k])
trough_r = np.array([r_at(k) for k in trough_k])

print("peaks:   ", list(zip(peak_k, peak_r)))
print("troughs: ", list(zip(trough_k, trough_r)))
\end{verbatim}

Part 2:

\begin{verbatim}
def fit_model(k_array, r_array, fix_A=None):
    """
    Fit r_k approx. A + B * (log log k / log k) + C * (1 / log k)
    on the given k_array, r_array (1D numpy arrays).

    If fix_A is not None, A is held fixed at that value and only B,C are fit.
    Returns (A, B, C, R2, r_pred).
    """
    k = np.asarray(k_array, dtype=float)
    y = np.asarray(r_array, dtype=float)

    x1 = np.log(np.log(k)) / np.log(k)   # term for B
    x2 = 1.0 / np.log(k)                 # term for C

    if fix_A is None:
        # full 3-parameter fit
        X = np.column_stack([np.ones_like(k), x1, x2])
        coeffs, residuals, rank, s = np.linalg.lstsq(X, y, rcond=None)
        A, B, C = coeffs
        y_pred = X @ coeffs
    else:
        # fix A, fit B,C
        A = float(fix_A)
        X = np.column_stack([x1, x2])
        y_shift = y - A
        coeffs, residuals, rank, s = np.linalg.lstsq(X, y_shift, rcond=None)
        B, C = coeffs
        y_pred = A + X @ coeffs

    # R^2
    y_mean = y.mean()
    ss_tot = np.sum((y - y_mean)**2)
    ss_res = np.sum((y - y_pred)**2)
    R2 = 1.0 - ss_res/ss_tot if ss_tot != 0 else 1.0

    return A, B, C, R2, y_pred

# fit peaks with A free
A_p, B_p, C_p, R2_p, pred_peak = fit_model(peak_k, peak_r, fix_A=None)
print("Peaks (A free): A =", A_p, "B =", B_p, "C =", C_p, "R^2 =", R2_p)

# fit peaks with A fixed at 2
A_p2, B_p2, C_p2, R2_p2, pred_peak2 = fit_model(peak_k, peak_r, fix_A=2.0)
print("Peaks (A=2):    B =", B_p2, "C =", C_p2, "R^2 =", R2_p2)

# fit troughs (A free)
A_t, B_t, C_t, R2_t, pred_trough = fit_model(trough_k, trough_r, fix_A=None)
print("Troughs: A =", A_t, "B =", B_t, "C =", C_t, "R^2 =", R2_t)

# fit troughs (A = 2)
A_t2, B_t2, C_t2, R2_t2, pred_trough = fit_model(trough_k, trough_r, fix_A=2.0)
print("Troughs (A=2): B =", B_t2, "C =", C_t2, "R^2 =", R2_t2)
\end{verbatim}

Output of part 2 (wrapped):

\begin{verbatim}
Peaks (A free): A = 1.9471086644613707 B = -0.6635108713154201 C = -0.9061856801140058
R^2 = 0.9991547825053999
Peaks (A=2):    B = -1.0503853115520978 C = -0.5288618840582724 R^2 = 0.99881703775137
Troughs: A = 1.8260681683738214 B = -0.43568639887707306 C = -0.7093619926349914
R^2 = 0.9980947124803357
Troughs (A=2): B = -1.729749295822425 C = 0.5773223100098374 R^2 = 0.9921324249943385
\end{verbatim}

Part 3:

\begin{verbatim}
# Delete last element

# fit peaks with A free
A_p, B_p, C_p, R2_p, pred_peak = fit_model(peak_k[:-1], peak_r[:-1], fix_A=None)
print("Peaks (A free): A =", A_p, "B =", B_p, "C =", C_p, "R^2 =", R2_p)

# fit peaks with A fixed at 2
A_p2, B_p2, C_p2, R2_p2, pred_peak2 = fit_model(peak_k[:-1], peak_r[:-1], fix_A=2.0)
print("Peaks (A=2):    B =", B_p2, "C =", C_p2, "R^2 =", R2_p2)

# fit troughs (A free)
A_t, B_t, C_t, R2_t, pred_trough = fit_model(trough_k[:-1], trough_r[:-1], fix_A=None)
print("Troughs: A =", A_t, "B =", B_t, "C =", C_t, "R^2 =", R2_t)

# fit troughs (A = 2)
A_t2, B_t2, C_t2, R2_t2, pred_trough = fit_model(trough_k[:-1], trough_r[:-1], fix_A=2.0)
print("Troughs (A=2): B =", B_t2, "C =", C_t2, "R^2 =", R2_t2)
\end{verbatim}

Output of part 3 (wrapped):

\begin{verbatim}
Peaks (A free): A = 1.98416458371712 B = -0.9175392149489754 C = -0.6735972057449692
R^2 = 0.9992440762971366
Peaks (A=2):    B = -1.0305224682459564 C = -0.565939141120064 R^2 = 0.9992203523314546
Troughs: A = 1.8159989832953467 B = -0.3652401650649837 C = -0.7753932406157733
R^2 = 0.9974343578780224
Troughs (A=2): B = -1.6982759403090761 C = 0.5177648430783589 R^2 = 0.9923355870025801
\end{verbatim}

Part 4:

\begin{verbatim}
# Delete first element

# fit peaks with A free
A_p, B_p, C_p, R2_p, pred_peak = fit_model(peak_k[1:], peak_r[1:], fix_A=None)
print("Peaks (A free): A =", A_p, "B =", B_p, "C =", C_p, "R^2 =", R2_p)

# fit peaks with A fixed at 2
A_p2, B_p2, C_p2, R2_p2, pred_peak2 = fit_model(peak_k[1:], peak_r[1:], fix_A=2.0)
print("Peaks (A=2):    B =", B_p2, "C =", C_p2, "R^2 =", R2_p2)

# fit troughs (A free)
A_t, B_t, C_t, R2_t, pred_trough = fit_model(trough_k[1:], trough_r[1:], fix_A=None)
print("Troughs: A =", A_t, "B =", B_t, "C =", C_t, "R^2 =", R2_t)

# fit troughs (A = 2)
A_t2, B_t2, C_t2, R2_t2, pred_trough = fit_model(trough_k[1:], trough_r[1:], fix_A=2.0)
print("Troughs (A=2): B =", B_t2, "C =", C_t2, "R^2 =", R2_t2)
\end{verbatim}

Output of part 4 (wrapped):

\begin{verbatim}
Peaks (A free): A = 1.899948394621106 B = -0.2899403839734418 C = -1.3024010728889466
R^2 = 0.9990412555574135
Peaks (A=2):    B = -1.0539592912295694 C = -0.5213961090436496 R^2 = 0.9981449829774532
Troughs: A = 1.8928451151056718 B = -0.9716702942910586 C = -0.13305801936408185
R^2 = 0.9978143593204098
Troughs (A=2): B = -1.804952644808632 C = 0.7353980066262816 R^2 = 0.9961298225748858
\end{verbatim}

\subsection*{Direct plot}

\begin{verbatim}
import numpy as np
import math
import matplotlib.pyplot as plt

# seq: numpy array [a_1, ..., a_N]
N = len(seq)

# we'll use k = 2,...,N (since log 1 = 0 and log log 1 is undefined)
ks = np.arange(2, N+1)
a_vals = seq[1:]   # a_2,...,a_N

# f(k) = log a_k - 2 log k + log log k
f_vals = np.log(a_vals) - 2.0 * np.log(ks) + np.log(np.log(ks))

plt.figure(figsize=(12,6))
plt.plot(np.log(ks)[20:], f_vals[20:], label=r"$\log a_k - 2\log k + \log\log k$")
plt.axhline(-0.64, color="black", linewidth=0.8, linestyle="--")
plt.axline((7, -0.84), (10, -1.14), color="black", linewidth=0.8, linestyle="--")
plt.xlabel(r"$\log k$")
plt.ylabel(r"$\log a_k - 2\log k + \log\log k$")
plt.title(r"Deviation from $a_k \sim k^2 / \log k$")
plt.legend()
plt.tight_layout()
plt.show()
\end{verbatim}

\end{document}